\documentclass{article}
\usepackage{amssymb}
\usepackage{amsmath}
\usepackage{amsfonts}
\usepackage{graphicx}
\usepackage[english]{babel}
\usepackage{graphics}

\newtheorem{theorem}{Theorem}[section]
\newtheorem{lemma}[theorem]{Lemma}
\newtheorem{corollary}[theorem]{Corollary}

\newtheorem{remark}[theorem]{Remark}

\newcommand{\qed}{\hfill \rule{2.3mm}{2.3mm}}

\newcommand{\al}{\alpha}

\newcommand{\ve}{\varepsilon}
\newcommand{\vp}{\varphi}

\newcommand*{\esssup}{\operatorname*{ess\phantom{|}\!sup}}
\newcommand*{\essinf}{\operatorname*{ess\phantom{|}\!inf}}

\newcommand{\bee}{\begin{equation}}
\newcommand{\ee}{\end{equation}}
\newcommand{\uu}{\bar{u}_\ve}

\newcommand{\R}{\mathbb{R}}

\newcommand{\Z}{\mathbb{Z}}

\newcommand{\x}{x/\varepsilon}

\newcommand{\reff}[1]{(\ref{#1})}

\title{\bf A common approach to singular perturbation and homogenization II:  Semilinear elliptic systems} 

\newcounter{thesame}
\author{
	Nikolai N. Nefedov
	\thanks{Lomonosov Moscow State University, Faculty of Physics, 19899 Moscow, Russia. 
		{\small   E-mail:
			{\tt nefedov@phys.msu.ru}
	}}
	\ \ \ Lutz Recke \thanks{Humboldt University of Berlin, Institute of Mathematics, Rudower Chausee 25, 12489 Berlin, Germany.
		{\small   E-mail:
			{\tt lutz.recke@hu-berlin.de}}
}}

\date{}

\begin{document}

\maketitle

\noindent
\begin{abstract}
\noindent
We consider  periodic homogenization of
boundary value problems for second-order semilinear elliptic systems in 2D of the type
$$
\partial_{x_i}\Big(a^{\al \beta}_{ij}(\x)
\partial_{x_j}u^\beta(x)
+b_i^\al(x,u(x))\Big)=b^\al(x,u(x))
\mbox{ for } x \in \Omega.
$$
For small $\ve>0$ we prove existence of weak solutions $u=u_\ve$ as well as their local uniqueness for $\|u-u_0\|_\infty \approx 0$,
where $u_0$ is a given non-degenerate
weak solution to the homogenized boundary value problem, and we estimate the rate of convergence to zero of
$\|u_\ve-u_0\|_\infty$
for $\ve \to 0$.

Our assumptions are, roughly speaking, as follows: The functions $a_{ij}^{\al \beta}$
are bounded, measurable and  
$\Z^2$-periodic, the functions 
$b_i^\al(\cdot,u)$ and
$b^\al(\cdot,u)$ are 
bounded and measurable, the functions
$b_i^\al(x,\cdot)$ and
$b^\al(x,\cdot)$
are $C^1$-smooth, and $\Omega$ is a bounded Lipschitz domain in $\R^2$.
Neither global solution uniqueness is supposed nor growth restrictions of $b_i^\al(x,\cdot)$
or
$b^\al(x,\cdot)$
nor higher regularity of $u_0$, and cross-diffusion is allowed.

The main tool of the proofs is an abstract result of implicit function theorem type
which in the past has been applied to singularly perturbed nonlinear ODEs and elliptic and parabolic PDEs and, hence, which permits a common approach to
existence, local uniqueness and error estimates for singularly perturbed problems and and for homogenization problems.

\end{abstract}  



\section{Introduction}
\label{sec:1}
\setcounter{equation}{0}
\setcounter{theorem}{0}

In this paper we present an abstract result of implicit function theorem type (see Section 2), which in the past has been applied to singularly pertubed nonlinear ODEs and PDEs in \cite{Butetc,But2022,Fiedler,NURS,
OmelchenkoRecke2015,Recke2022,
ReckeOmelchenko2008}
and, in Part I \cite{1}, to periodic homogenization of quasilinear ODE systems. In the present paper we apply it to periodic homogenization 
of Dirichlet problems for 2D semilinear elliptic PDE systems of the type
\bee
\label{ODE}
\left.
\begin{array}{l}
\partial_{x_i}\Big(a^{\al \beta}_{ij}(\x)
\partial_{x_j}u^\beta(x)
+b_i^\al(x,u(x))\Big)=b^\al(x,u(x)) \mbox{ for } x \in \Omega,\\
u^\al(x)=0\mbox{ for } x \in \partial\Omega,
\end{array}
\right\}
\al=1,\ldots,n
\ee
as well as of other boundary value problems for those systems (see Section \ref{sec4}).
Here and in the following repeated indices are to be summed over $\al,\beta,\gamma,\ldots=1,\ldots,n$ and $i,j,k,\ldots=1,2$,
and $\ve>0$ is the small homogenization parameter.
We assume that
\begin{eqnarray}
\label{Omass}
&&\mbox{$\Omega$ is a bounded Lipschitz domain in $\R^2$,}\\
\label{perass}
&&a_{ij}^{\al \beta} \in L^\infty(\R^2) \mbox{ and } a_{ij}^{\al \beta}(\cdot+z)=
a_{ij}^{\al \beta}\mbox{ for all } z \in \Z^2,\\
\label{diffass}
&&u\in \R^n \mapsto (b_i^\al(\cdot,u),
b^\al(\cdot,u))
\in L^\infty(\Omega;\R^2) \mbox{ is $C^1$-smooth}
\end{eqnarray}
and
\bee
\label{alass}
\inf\left\{\int_{\R^2}a_{ij}^{\al \beta}(y)
\partial_{y_i}\vp^\al(y)\partial_{y_j}\vp^\beta(y)dy:\;
\vp \in C_c^\infty(\R^2;\R^n),\;
\int_{\R^2}\partial_{y_i}\vp^\al(y)\partial_{y_i}\vp^\al(y)dy=1\right\}>0.
\ee
The components of the homogenized diffusion tensor are, by definition,
\bee
\label{Anullcoef}
\hat{a}^{\al\beta}_{ij}:=
\int_{(0,1)^2}\left(a^{\al \beta}_{ij}(y)
+a^{\al \gamma}_{ik}(y)
\partial_{y_k}v^{\gamma \beta}_j(y)\right)dy,
\ee
where the $2n$ correctors $v^{\beta}_j
\in W^{1,2}_{\rm loc}(\R^2;\R^n)$, $\beta=1,\ldots,n$, $j=1,2$,
are defined by the $2n$ cell problems
\bee
\label{cell}
\left.
\begin{array}{l}
\partial_{y_i}
\left(a^{\al \beta}_{ij}(y)+
a^{\al \gamma}_{ik}(y)
\partial_{y_k}v^{\gamma \beta}_j(y)\right)
=0\mbox{ for } y\in \R^2,\\
v^{\al \beta}_j(\cdot+z)=v^{\al \beta}_j \mbox{ for } z \in \Z^2,\;
\displaystyle\int_{(0,1)^2} v^{\al \beta}_j(y))dy=0,
\end{array}
\right\}
\al=1,\ldots,n.
\ee
It is well-known (as a consequence of 
assumption \reff{alass} and
the Lax-Milgram lemma, cf., e.g. \cite[Section 2.2 and Lemma 2.2.4]{Shen}) that  the problems \reff{cell} are uniquely weakly solvable  
and that the homogenized diffusion coefficients $\hat{a}^{\al\beta}_{ij}$ satisfy the coercivity condition
\reff{alass} as well, i.e.
\bee
\label{alass1}
\inf\left\{\int_{\R^2}\hat a_{ij}^{\al \beta}
\partial_{y_i}\vp^\al(y)\partial_{y_j}\vp^\beta(y)dy:\;
\vp \in C_c^\infty(\R^2;\R^n),\;
\int_{\R^2}\partial_{y_i}\vp^\al(y)\partial_{y_i}\vp^\al(y)dy=1\right\}>0.
\ee

Let us formulate our main result. It concerns existence and local uniqueness of
weak solutions $u=u_\ve$ to \reff{ODE} with $\ve \approx 0$, which are close to a given non-degenerate weak solution $u=u_0$ to
the homogenized problem
\bee
\label{hombvp}
\left.
\begin{array}{l}
\partial_{x_j}\Big(\hat{a}^{\al \beta}_{ij}
\partial_{x_j}u^\beta(x)
+b_i^\al(x,u(x))\Big)
=b^\al(x,u(x)) \mbox{ for } x \in \Omega,\\
u^\al(x)=0 \mbox{ for } x \in \partial\Omega,
\end{array}
\right\}
\al=1,\ldots,n
\ee
as well as the rate of convergence to zero for $\ve \to 0$ of the homogenization error $\|u_\ve-u_0\|_\infty$.
Here and in what follows we denote 
\bee
\label{infdef}
\|u\|_\infty:=\max_{\al=1,\ldots,n}\esssup_{x \in \Omega}|u^\al(x)|
\mbox{ for }
u \in L^\infty(\Omega;\R^n).
\ee
As usual, a vector function $u \in W_0^{1,2}(\Omega;\R^n)\cap C(\overline\Omega;\R^n)$ is called weak solution to the boundary value problem \reff{ODE} if it satisfies the variational equation
$$
\int_\Omega
\left(\Big(a^{\al \beta}_{ij}(\x)
\partial_{x_j}u^\beta(x)+b_i^\al(x,u(x))\Big)
\partial_{x_i}\vp^\al(x)+b^\al(x,u(x))\vp^\al(x)\right)dx=0
$$
for all 
$\vp \in W_0^{1,2}(\Omega;\R^n)$,
and similarly for \reff{hombvp} and for its linearization
\reff{linbvp} and for the cell problems \reff{cell}.

\begin{theorem} 
\label{main}
Suppose \reff{Omass}-\reff{alass}, and
let $u=u_0$ be a weak solution to \reff{hombvp} such that
the linearized homogenized boundary value problem
\bee
\label{linbvp}
\left.
\begin{array}{l}
\partial_{x_j}\Big(\hat{a}^{\al \beta}_{ij}
\partial_{x_j}u^\beta(x)
+\partial_{u^\gamma}b_i^\al(x,u_0(x))u^\gamma(x)\Big)
=\partial_{u^\gamma}b^\al(x,u_0(x))u^\gamma(x) \mbox{ for } x \in \Omega,\\
u^\al(x)=0 \mbox{ for } x \in \partial\Omega,
\end{array}
\right\}
\ee
for $\al=1,\ldots,n$ does not have weak solutions $u\not=0$.

Then the following is true:

(i)
There exist $\ve_0>0$ and  $\delta>0$
such that for all $\ve \in (0,\ve_0]$ there exists exactly one 
weak solution $u=u_\ve$ 
to \reff{ODE} with $\|u-u_0\|_\infty \le \delta$. Moreover, 
$$
\|u_\ve-u_0\|_\infty\to 0 \mbox{ for }
\ve \to 0.
$$

(ii) If 
$u_0\in W^{2,p_0}(\Omega;\R^n)$ with certain $p_0>2$, then for any $p>2$ we have
\bee
\label{urate}
\|u_\ve-u_0\|_\infty=O(\ve^{1/p}) 
\mbox{ for } \ve \to 0.
\ee
\end{theorem}

\begin{remark}
\label{unullreg}
Weak solutions to \reff{hombvp}
with \reff{Omass}-\reff{alass} 
do not belong to 
$W^{2,p_0}(\Omega;\R^n)$ with certain $p_0>2$,
in general. But they do, if the boundary $\partial \Omega$ is $C^{1,1}$-smooth (or if $\Omega$ is convex), 
and if the drift functions $b_i^\al$ are $C^{1}$-smooth.
This problem is well-known in homogenization theory (cf. e.g.
\cite[Remark 6.4]{Ci}, \cite[Remark 4.2]{Senik}, \cite[Section 3]{QXu}).
\end{remark}

\begin{remark}
\label{Campa}
We conjecture that Theorem \ref{main} remains to be true for any space dimension if the elliptic system is triangular, i.e. $a^{\al \beta}_{ij}=0$ for $\al>\beta$ (in particular, for scalar elliptic equations), this is work in progress \cite{GRIII}.
The reason for that is the following:
In the present paper, which concerns space dimension two, we use K. Gr\"ogers result \cite{G} (cf. Theorem \ref{Gro} of the present paper)
about   maximal regularity of boundary value problems for elliptic systems
with non-smooth data in the pair of Sobolev spaces $W_0^{1,p}(\Omega;\R^n)$ and $W^{-1,p}(\Omega;\R^n)$ with $p\approx 2$
as well as the continuous embedding $W^{1,p}(\Omega;\R^n)\hookrightarrow C(\overline\Omega;\R^n)$ for  $p>2$
in the case of space dimension two. But there exists a replacement of these results for triangular systems with any space dimension
(cf. \cite{GR,GrR}), where the Sobolev spaces are replaced by appropriate Sobolev-Campanato spaces which for any space dimension are continuously embedded into $C(\overline\Omega;\R^n)$.

Remark that in \cite{GR,G,GrR} 
more general types of boundary conditions are allowed, for example mixed Dirichlet-Robin boundary conditions. Hence,
Theorem \ref{main} is true also for those boundary conditions (cf. Section  \ref{sec4} of the present paper), and
we expect the same in the case of any 
space dimension and $\hat a^{\al \beta}_{ij}=0$ for $\al>\beta$.

What concerns elliptic systems of the type \reff{ODE} with space dimension larger than two, which are far from being triangular, we do not believe that results of the type of Theorem \ref{main} are true, in general.
\end{remark}

\begin{remark}
\label{sufficient}
The assumption of Theorem \ref{main}, that there do not exist nontrivial weak solutions to \reff{hombvp}, is rather implicit. But there exist simple explicit sufficient conditions for it. For example, if 
$$
\inf\left\{\int_{\Omega} \partial_{u^\beta}b^\al(x,u_0(x))u^\al(x)u^\beta(x)dx:\;
u \in W_0^{1,2}(\Omega;\R^n),\;
\int_{\Omega} \partial_{x_i}u^\al(x)\partial_{x_i}u^\al(x)dx
=1\right\}>0,
$$
and if this infimum as well as the infimum in \reff{alass1} (which equals to the 
infimum in \reff{alass}, cf. \cite[Lemma 2.2.4]{Shen})
are sufficiently large in comparison with
$$
\sup\left\{\left|\int_{\Omega} \partial_{u^\beta}b_i^\al(x,u_0(x))\partial_{x_i}u^\al(x)u^\beta(x)dx\right|:\;
u \in W_0^{1,2}(\Omega;\R^n),\;
\int_{\Omega} \partial_{x_i}u^\al(x)\partial_{x_i}u^\al(x)dx
=1\right\},
$$
then there do not exist nontrivial weak solutions to \reff{hombvp}.
\end{remark}

\begin{remark}
\label{bprop}
In many applications the reaction functions $b^\al$ are of the type
$$
b^\al(x,u)=\sum_{l=1}^mc_l^\al(x)d_l^\al(u)
\mbox{ with } c_l^\al \in L^\infty(\Omega),\; d_l^\al \in C^1(\R^n),
$$
and those satisfy assumption \reff{diffass},
and similarly for the drift functions $b_i^\al$.
\end{remark}

\begin{remark}
\label{Aprop}
The integral condition \reff{alass} is often referred as $V$-ellipticity or $V$-coercivity, and it follows from the Legendre condition
\bee
\label{Leg}
\essinf\left\{a_{ij}^{\al \beta}(y)\xi_i^\al\xi_j^\beta:\; y\in \R^2,\;
\xi \in \R^{2n},\; \xi_i^\al\xi_i^\alpha=1\right\}>0
\ee
and it
implies the Legendre-Hadamard condition
$$
\essinf\left\{
a_{ij}^{\al \beta}(y)\xi_i\xi_j\eta^\al\eta^\beta:\;
y\in \R^2,\;
\xi\in \R^2,\;\eta \in \R^{n},\; \xi_i\xi_i=\eta^\al\eta^\al=1\right\}>0.
$$
If $n=1$ (i.e. in the case of scalar equations) or if the coefficients $a_{ij}^{\al \beta}$ are constant (as in \reff{hombvp}), then $V$-ellipticity is equivalent to the
Legendre-Hadamard condition.
\end{remark}

\begin{remark}
\label{Neuk}
For $L^\infty$  estimates of $u_\ve-u_0$ for 
periodic homogenization of
scalar linear elliptic PDEs 
see, e.g. \cite[Chapter 2.4]{Ben} and \cite{He} and for linear elliptic systems
\cite[Theorem 3.4]{Ke},
 \cite[Theorem 7.5.1]{Shen} and
\cite[Theorem 1.7]{WangZhang}.
For $L^p$ homogenization error estimates for scalar linear elliptic PDEs 
see, e.g. \cite{He} and \cite[Theorem 1.1]{QXu} and for linear elliptic systems \cite[Theorem 7.5.1]{Shen}
and \cite[Theorem 1.5]{QXu1}.

What concerns existence and local uniqueness
for nonlinear elliptic homogenization problems
(without assumption of global uniqueness) we know only the result 
\cite{Bun}  for scalar semilinear elliptic PDEs of the type
$\mbox{\rm div}\, a(\x) \nabla u(x)=f(x)g(u(x))$,
where the nonlinearity $g$ is supposed to have a sufficiently small local Lipschitz constant (on an appropriate bounded interval). Let us mention also \cite{Lanza1,Lanza2}, where existence and local uniqueness for a homogenization problem for the linear Poisson equation with nonlinear Robin boundary conditions
in a periodically perforated domain is shown. There the specific structure of the problem (no highly oscillating diffusion coefficients) allows to apply the classical implicit function theorem.

For periodic homogenization of linear elliptic PDEs (with small homogenization parameter $\ve$), which are singularly perturbed (with small singular perturbation parameter $\delta$) see \cite{Tev}.
\end{remark}

~\\ 

Our paper is organized as follows: 
In Section \ref{secabstract} we consider abstract nonlinear parameter depending equations of the type
\bee
\label{intrabstract}
F_\ve(u)=0.
\ee
Here $\ve>0$ is the parameter. We prove a result on existence and local uniqueness of a family of solutions $u=u_\ve \approx \bar u_\ve$ to \reff{intrabstract} with $\ve \approx 0$, where $\bar u_\ve$ is a family of  approximate solutions to \reff{intrabstract}, i.e. a family with
$F_\ve(\bar u_\ve)\to 0$ for $\ve \to 0$, and we estimate the norm of the error $u_\ve-\bar u_\ve$ by the norm of the discrepancy $F_\ve(\bar u_\ve)$. 
This type of generalized implicit function theorems has been successfully applied to singularly perturbed nonlinear ODEs and PDEs and to homogenization of nonlinear ODEs.
Contrary to the classical implicit function theorem it is not supposed that
the linearized operators $F'_\ve(u)$ converge for $\ve \to 0$ 
with respect to  the uniform operator norm. And, indeed, in the 
applications to singularly perturbed problems as well as to  homogenization 
problems they do not converge for $\ve \to 0$ 
with respect to  the uniform operator norm (cf. Remark \ref{linop} below).
Hence, the present paper introduces an application 
to semilinear elliptic PDE systems of a common approach to  existence, local uniqueness and error estimates for singularly perturbed problems and for homogenization problems.
Another application to periodic homogenization of quasilinear ODE systems
of the type
$$
a(x,\x,u(x),u'(x))'=f(x,\x,u(x))
$$
has been presented in Part I \cite{1}.

In Section \ref{sec3} we prove Theorem \ref{main} by means of 
the results of Section~\ref{secabstract}.
Here the main work is to construct
appropriate families of approximate  solutions 
to \reff{ODE} with $\ve \approx 0$
with small discrepancies in appropriate function space norms.
In order to apply implicit function theorems mainly one needs isomorphism properties 
of the linearized operators. In the setting of Section~\ref{sec3} they follow from
K. Gr\"ogers result \cite{G}
about   maximal regularity of boundary value problems for elliptic systems
with non-smooth data in the pair of Sobolev spaces $W_0^{1,p}(\Omega;\R^n)$ and $W^{-1,p}(\Omega;\R^n)$ with $p\approx 2$.
In order to apply implicit function theorems one needs also $C^1$-smoothness of the appearing nonlinear superposition operators. In the
setting of Section~\ref{sec3} these operators have to be well-defined and $C^1$-smooth on the Sobolev spaces $W^{1,p}(\Omega;\R^n)$
with $p>2$, but $p\approx 2$, and therefore we have to suppose that the space dimension is two.

Finally, in Section \ref{sec4} we show that Theorem \ref{main} is true also in the case of mixed boundary conditions. More exactly, we consider homogeneous Dirichlet boundary conditions on one part of $\partial \Omega$ and inhomogeneous  natural (i.e. nonlinear Robin) boundary conditions on the remaining part, and the two boundary parts are allowed to touch.

\section{An abstract result of implicit function theorem type}
\label{secabstract}
\setcounter{equation}{0}
\setcounter{theorem}{0}

Let $U$ and $V$ be a Banach spaces with norms $\|\cdot\|_U$ and $\|\cdot\|_V$, respectively.
For $\ve>0$ let be given 
$$
\uu \in U \mbox{ and } F_\ve\in C^1(U;V).
$$
In this section we consider the abstract equation
\bee
\label{abeq}
F_\ve(u)=0. 
\ee
Roughly speaking, we will show the following:
If the elements $\uu$ satisfy \reff{abeq} approximately for $\ve \approx 0$, i.e. if $\|F_\ve(\uu)\|_V\to 0$ for  $\ve \to 0$, 
and if they are non-degenerate solutions
(cf. assumption \reff{coerz} below), then
for $\ve \approx 0$ there exists exactly one solution $u=u_\ve$ to \reff{abeq} with $\|u-\uu\|_U\approx 0$, and $\|u_\ve-\uu\|_U=O(\|F_\ve(\uu)\|_V)$ for $\ve \to 0$.
For that we do not suppose any convergence of the  operators $F_\ve$ or $F'_\ve(u)$ or of the elements $\uu$ for $\ve \to 0$.
Remark that in the classical implicit function theorem one cannot omit, in general, the assumption, that $F_\ve'(u)$ converges for $\ve \to 0$ with respect to the uniform operator norm (cf. \cite[Section 3.6]{Katz}).

\begin{theorem}
\label{ift}
Suppose that 
\bee
\label{coerz}
\left.
\begin{array}{l}
\mbox{there exist $\ve_0>0$ and $c>0$ such that for all $\ve \in (0,\ve_0]$ the operators $F_\ve'(\uu)$ are}\\
\mbox{Fredholm of index zero from }
$U$ \mbox{ into } $V$,
\mbox{and } \|F_\ve'(\uu))u\|_V \ge c \|u\|_U \mbox{ for all } u \in U
\end{array}
\right\}
\ee
and
\bee
\label{Fas}
\sup_{\|v\|_U \le 1}\|(F'_\ve(\uu+u)-F'_\ve(\uu))v\|_V \to 0
\mbox{ for } \ve + \|u\|_U \to 0
\ee
and
\bee
\label{lim}
\|F_\ve(\uu)\|_V \to 0 \mbox{ for } \ve \to 0.
\ee

Then there exist $\ve_1 \in (0,\ve_0]$ and $\delta>0$  such that for all $\ve \in (0,\ve_1]$ 
there exists exactly one  solution $u = u_\ve$ to~\reff{abeq} with $\|u_\ve-\uu\|_U \le \delta$, and 
\bee
\label{apriori1}
\|u_\ve-\uu\|_U \le \frac{2}{c}\|F_\ve(\uu))\|_V.
\ee
\end{theorem}
{\bf Proof }
Take $\ve \in (0,\ve_0]$. 
Because of assumption \reff{coerz} the linear operator $F'_\ve(\uu)$ is an isomorphism from $U$ onto $V$, and, hence,
equation \reff{abeq} is equivalent to the fixed point problem
$$
u=G_\ve(u):=u-F_\ve'(\uu))^{-1}F_\ve(u).
$$
Take $u_1,u_2 \in U$. Then 
\begin{eqnarray}
\label{strict}
\|G_{\ve}(u_1) - G_{\ve}(u_2)\|_U
&=&\left\|F_\ve'(\uu)^{-1}\int_0^1\left(
F'_\ve(\uu)-
F'_{\ve}(s u_1 + (1-s) u_2)\right)ds (u_1 - u_2)\right\|_U\nonumber\\
&\le& \frac{1}{c}\, \max_{0 \le s \le 1}\|\left(F'_\ve(\uu)-F'_{\ve}(s u_1 + (1-s) u_2)\right) (u_1 - u_2)\|_V.
\end{eqnarray}  
Here we used that \reff{coerz} yields that
$c\|F_\ve'(\uu)^{-1}v\|_U \le \|v\|_V$
for all $v \in V$.

Denote ${\cal B}^r_\ve:=\{u \in U:\; \|u-\uu\|_U \le r\}$.
If $u_1,u_2 \in {\cal B}^r_\ve$,
then also $su_1+(1-s)u_2 \in {\cal B}^r_\ve$
for all $s \in [0,1]$. Therefore  it follows from \reff{Fas} and \reff{strict} that there exist $r_0>0$ and $\ve_1 \in (0,\ve_0]$ such that for all $\ve \in (0,\ve_1]$
the maps $G_\ve$ are strictly contractive with contraction constant $1/2$ on the closed balls ${\cal B}^{r_0}_\ve$.
Moreover, for all $\ve \in (0,\ve_1]$ and $u \in {\cal B}^{r_0}_\ve$ we have
\begin{eqnarray*}
\label{in}
\left\|G_{\ve}(u) - \uu\right\|_U &\le & \left\|G_{\ve}(u) - G_\ve(\uu)\right\|_U
+\left\|G_{\ve}(\uu) - \uu\right\|_U\\
&\le & \frac{r_0}{2}+\left\|F_\ve'(\uu))^{-1}F_\ve(\uu)\right\|_U
\le \frac{r_0}{2}+\frac{1}{c}\left\|F_\ve(\uu)\right\|_V,
\end{eqnarray*}
and \reff{lim} yields that $G_\ve$ maps ${\cal B}^{r_0}_\ve$ into itself if $\ve_1$ is taken sufficiently small.

Now, Banach's fixed point principle yields the existence and uniqueness assertions of Theorem \ref{ift}, and the estimate \reff{apriori1}
follows as above:
$$
\|u_\ve-\uu\|_U \le \|G_{\ve}(u_\ve) - G_\ve(\uu)\|_V+\|G_{\ve}(\uu) - \uu\|_U
\le \frac{1}{2}\|u_\ve-\uu\|_V
+\frac{1}{c}\left\|F_\ve(\uu)\right\|_V.
$$
\qed

\begin{remark}
In \cite{Butetc,But2022,Fiedler,NURS,
OmelchenkoRecke2015,Recke2022,
ReckeOmelchenko2008}) various  versions of Theorem \ref{ift} are presented. They differ slightly according to which problems they are applied (ODEs or elliptic or parabolic PDEs, stationary or time-periodic solutions, semilinear or quasilinear problems, smooth or nonsmooth data, one- or multi-dimensional perturbation parameter $\ve$).

For another result of the type of Theorem \ref{ift} and its applications to semilinear elliptic PDE systems with numerically  determined approximate solutions  see \cite[Theorem 2.1]{Breden}.
\end{remark}

If one applies Theorem \ref{ift}, for example to boundary value problems for elliptic PDEs, then different choices of function spaces $U$ and $V$ and of their norms $\|\cdot\|_U$ and $\|\cdot\|_V$ 
and of the family $\uu$ of approximate solutions
are appropriate. Criteria for these choices often are the following: The family $\uu$ of should be "simple" (for example, $\uu$ should be $\ve$-independent or given more less explicit in closed formulas, or to determine $ \uu$ numerically should be much cheeper than 
to determine the exact solution $u_\ve$ numerically), and the rate of convergence to zero of $\|F_\ve(\uu)\|_V$ for $\ve \to 0$ should be high. The norm $\|\cdot\|_V$ should be weak and the norm $\|\cdot\|_U$ should be strong such that the error estimate 
\reff{apriori1}
is strong. But at the same time 
the norm $\|\cdot\|_U$ should be weak such that the domain of local uniqueness, which contains all $u \in U$ with $\|u-\uu\|_U\le \delta$, is large.
These criteria are contradicting, of course. Hence, in any application of Theorem \ref{ift}
the choices of $U$, $V$, $\|\cdot\|_U$, $\|\cdot\|_V$ and $\uu$ are compromises according to the requirements of the application.

One way to find such compromises is described in  Corollary \ref{cor} below. It delivers existence and local uniqueness of solutions $u=u_\ve$ to the equation $F_\ve(u)=0$ with $\ve \approx 0$ and $\|u-u_0\|\approx 0$, where $\|\cdot\|$ is another norm in $U$, which is allowed to be much weaker
than the norm $\|\cdot\|_U$,
and
where $\|F_\ve(u_0)\|_V$ does not converge to zero for $\ve \to 0$, in general. The price for that is that the estimate \reff{apriori2} of the error $u_\ve-u_0$ is with respect to the weaker norm $\|\cdot\|$, only.

\begin{corollary}
\label{cor}
Suppose
\reff{lim}. Further, let be given $u_0 \in U$ 
such that
\bee
\label{coerza}
\left.
\begin{array}{l}
\mbox{there exist $\ve_0>0$ and $c>0$ such that for all $\ve \in (0,\ve_0]$ the operators $F_\ve'(u_0)$ are}\\
\mbox{Fredholm of index zero from }
$U$ \mbox{ into } $V$,
\mbox{and } \|F_\ve'(u_0))u\|_V \ge c \|u\|_U \mbox{ for all } u \in U,
\end{array}
\right\}
\ee
and let be given a norm $\|\cdot\|$ in $U$ 
such that
\begin{eqnarray}
\label{weaker}
&& d:=\sup\{\|u\|: \; u \in U, \|u\|_U \le 1\}< \infty,\\
\label{newconv}
&& \|\uu-u_0\| \to 0 \mbox{ for } 
\ve \to 0,\\
\label{Fas1}
&&\sup_{\|v\|_U \le 1}\|(F'_\ve(u_0+u)-F'_\ve(u_0))v\|_V \to 0
\mbox{ for } \ve + \|u\| \to 0.
\end{eqnarray}

Then
there exist $\ve_1 \in (0,\ve_0]$ and $\delta>0$  such that for all $\ve \in (0,\ve_1]$ 
there exists exactly one  solution $u = u_\ve$ to~\reff{abeq} with $\|u-u_0\| \le \delta$,
and
\bee
\label{apriori2}
\|u_\ve-u_0\| \le \|\uu-u_0\|+
\frac{3d}{c}\|F_\ve(\uu))\|_V.
\ee
\end{corollary}
{\bf Proof }
For all $\ve>0$ and $u,v \in U$ we have
$$
\|(F_\ve'(\uu+u)-F_\ve'(\uu))v\|_V
\le
\|(F_\ve'(u_0+(\uu-u_0)+u)-F_\ve'(u_0))v\|_V+
\|(F_\ve'(u_0)-F_\ve'(u_0+(\uu-u_0)))v\|_V.
$$
Hence, assumptions \reff{weaker}-\reff{Fas1} imply that \reff{Fas} is satisfied.
Similarly, for all $\ve>0$ and $v \in U$ we have
$$
\|(F_\ve'(\uu)v\|_V
\ge
\|(F_\ve'(u_0)v\|_V-
\|(F_\ve'(u_0)-F_\ve'(u_0+(\uu-u_0)))v\|_V.
$$
Therefore assumptions \reff{coerza}, \reff{newconv} and \reff{Fas1} 
imply that \reff{coerz} is satisfied
(with another $\ve_0$ in \reff{coerz} than in \reff{coerza} and with arbitrary smaller $c$ in \reff{coerz}
than in \reff{coerza}).
Hence, Theorem \ref{ift} yields the existence assertion of Corollary \ref{cor} and the error estimate
$$
\|u_\ve-u_0\| \le \|\uu-u_0\|+d\|u_\ve-\uu\|_U\le \|\uu-u_0\|+\frac{3d}{c}\|F_\ve(\uu)\|_V.
$$

Now let us prove the local uniqueness assertion of Corollary \ref{cor}.
Take $\ve \in (0,\ve_1]$ and a solution $u \in U$ to \reff{abeq}. Then
$$
0=F_\ve(u)=F_\ve(\uu)+F'_\ve(\uu)(u-\uu)+
\int_0^1\left(F'_\ve(su+(1-s)\uu)-F'_\ve(\uu)\right)(u-\uu)ds,
$$
i.e.
$$
u-\uu=-F'_\ve(\uu)^{-1}\left(
F_\ve(\uu)+\int_0^1\left(F'_\ve(su+(1-s)\uu)-F'_\ve(\uu)\right)(u-\uu)ds\right),
$$
i.e.
\bee
\label{unest}
\|u-\uu\|_U\le \frac{1}{c}\left(
\|F_\ve(\uu)\|_V+\max_{0\le s \le 1}\|\left(F'_\ve(su+(1-s)\uu)-F'_\ve(\uu)\right)(u-\uu)\|_V\right).
\ee
But \reff{Fas1} yields that
$$
\max_{0\le s \le 1}\|\left(F'_\ve(su+(1-s)\uu)-F'_\ve(\uu)\right)(u-\uu)\|_V
=o(\|u-\uu\|_U) \mbox{ for }
\ve +\|u-\uu\| \to 0.
$$
Therefore \reff{newconv} implies that
$$
\max_{0\le s \le 1}\|\left(F'_\ve(su+(1-s)\uu)-F'_\ve(\uu)\right)(u-\uu)\|_V
=o(\|u-\uu\|_U) \mbox{ for }
\ve +\|u-u_0\| \to 0.
$$
Hence, if $\ve$ and $\|u-u_0\|$ are sufficiently small, then \reff{unest} yields that $\|u-\uu\|_U$ is sufficiently small, and the local uniqueness assertion of Theorem \reff{ift} implies that $u=u_\ve$.
\qed
\begin{remark}
In most of the applications of Corollary \ref{cor} to PDEs the element $u_0$ and the norm $\|\cdot\|$ are a priori given
(because for
$\ve \approx 0$ 
one is looking for solutions $u$ with
$\|u-u_0\|_\infty \approx 0$), and one has to choose  Banach spaces $U$ and $V$ (with their norms $\|\cdot\|_U$ and $\|\cdot\|_V$)
such that the PDE problem is equivalent to an abstract equation $F_\ve(u)=0$ with  $F_\ve
\in C^1(U;V)$
and with \reff{coerza}, \reff{weaker} and \reff{Fas1},
and one has to construct a family $\bar u_\ve$ with the properties \reff{lim} and \reff{newconv}.
But at the beginning one does not know if existence and local uniqueness for $\ve \approx 0$ and $\|u-u_0\|\approx 0$ is true or not for the given PDE problem. If not, then one is trying to choose and to construct something, which does not exist.

For example, in Theorem \ref{main}, which is the result of an application of Corollary \ref{cor} to problem \reff{ODE}, the spaces $U$ and $V$ and the family $\bar u_\ve$ are hidden only, they do not appear in the formulation of Theorem \ref{main}.
\end{remark}

\section{Proof of Theorem  \ref{main}}
\label{sec3}
\setcounter{equation}{0}
\setcounter{theorem}{0}
In this section we will prove Theorem \ref{main}
by means of Corollary \ref{cor}. 
For that we use the objects of Theorem \ref{main}: The bounded Lipschitz domain $\Omega \subset \R^2$, the diffusion coefficients $a^{\al \beta}_{ij}\in L^\infty(\R^2)$ with \reff{perass} and \reff{alass}, the drift and reaction functions $b_i^{\al},b^\al:\Omega\times \R^n\to \R$ with
\reff{diffass}, 
the correctors $v^{\beta}_j \in W^{1,2}_{\rm loc}(\R^2;\R^n)$, which are defined by the cell problems \reff{cell},
the homogenized diffusion coefficients
$\hat a^{\al \beta}_{ij} \in \R$, which are
defined in \reff{Anullcoef} and which satisfy \reff{alass1}, and
 the
weak solution $u_0\in W_0^{1,2}(\Omega;\R^n)
\cap C(\overline\Omega;\R^n)$ to the homogenized boundary value problem
\reff{hombvp}.

As usual, the norm in the Sobolev space $W^{1,p}(\Omega;\R^n)$ (with $p\ge 2$) is denoted by 
\bee
\label{1pnorm}
\|u\|_{1,p}:=
\left(\sum_{\al=1}^n
\int_\Omega
\left(|u^\al(x)|^p+
\sum_{i=1}^2 
|\partial_{x_i}u^\al(x)|^p\right)
dx\right)^{1/p},
\ee
$W_0^{1,p}(\Omega;\R^n)$ is the closure with respect to this norm 
of the set 
of all $C^{\infty}$-maps $u:\Omega \to \R^n$ with compact support,
and $W^{-1,p}(\Omega;\R^n):=W_0^{1,q}(\Omega;\R^n)^*$ is the dual space to 
$W_0^{1,q}(\Omega;\R^n)$ with $1/p+1/q=1$
with dual space norm
$$
\|f\|_{-1,p}:=\sup\{\langle f,\vp\rangle_{1,q}:\; \vp \in W_0^{1,q}(\Omega;\R^n),
\|\vp\|_{1,q}\le1\},
$$
where
$
\langle \cdot,\cdot\rangle_{1,q}:
W^{-1,p}(\Omega;\R^n)\times
W_0^{1,q}(\Omega;\R^n)\to \R
$
is the dual pairing. 

Further, we introduce linear bounded operators $A_0:W^{1,2}(\Omega;\R^n)\to W^{-1,2}(\Omega;\R^n)$ and, for $\ve>0$,
$A_\ve:W_0^{1,2}(\Omega;\R^n)\to W^{-1,2}(\Omega;\R^n)$
by
\bee
\label{Avedef}
\left.
\begin{array}{l}
\displaystyle\langle A_0 u,\vp\rangle_{1,2}:=
\int_\Omega \hat a_{ij}^{\al \beta}
\partial_{x_j}u^\beta(x)\partial_{x_i}\vp^\al(x)dx,\\
\displaystyle\langle A_\ve u,\vp\rangle_{1,2}:=
\int_\Omega a_{ij}^{\al \beta}(\x)
\partial_{x_j}u^\beta(x)\partial_{x_i}\vp^\al(x)dx,
\end{array}
\right\}
\mbox{ for all } \vp \in W_0^{1,2}(\Omega;\R^n).
\ee
Because of assumption \reff{perass} and the H\"older inequality we have the following:
For any $p \ge 2$ the restrictions of $A_0$ and $A_\ve$ to
$W^{1,p}(\Omega;\R^n)$
map $W^{1,p}(\Omega;\R^n)$ into
$W^{-1,p}(\Omega;\R^n)$, and
$$
\|A_0 u\|_{-1,p}+
\|A_\ve u\|_{-1,p}\le \mbox{const}\;\|u\|_{1,p} \mbox{ for all } u \in W_0^{1,p}(\Omega;\R^n),
$$
where the constant does not depend on $\ve$, $p$ and $u$. 
Further, because of assumption \reff{alass} we have that 
$$
\inf\left\{\langle A_\ve u,u\rangle_{1,2}:\;
\ve>0,\;
u \in W_0^{1,2}(\Omega;\R^n),
\; \|u\|_{1,2}=1\right\}>0,
$$
and similarly for $A_0$.
Therefore K. Gr\"oger's maximal regularity results for elliptic systems with non-smooth data \cite[Theorems 1 and 2 and Remark 14]{G} imply the following:
\begin{theorem}
\label{Gro}
There exist $p_1>2$ and $c>0$ such that
for all $\ve>0$ and all $p\in[2,p_1]$
the linear operators $A_0$ and $A_\ve$
are bijective from $W_0^{1,p}(\Omega;\R^n)$ onto $W^{-1,p}(\Omega;\R^n)$ and that
\bee
\label{Ainvert}
\|A_0^{-1}f\|_{1,p}+ 
\|A_\ve^{-1}f\|_{1,p} \le c\, \|f\|_{-1,p}
\mbox{ for all } f \in W^{-1,p}(\Omega;\R^n).
\ee
\end{theorem}

\begin{remark}
\label{linop}
(i) Estimates of the type \reff{Ainvert} often are called Meyers' estimates
(or estimates of Gr\"oger-Meyers type) because of the initiating paper \cite{M}
of N.G. Meyers. For the case of smooth boundaries $\partial \Omega$ see 
\cite[Theorem 4.1]{Ben}. See also \cite{CP} for $\Omega$ being a cube and for continuous diffusion coefficients as well as for certain transmission problems and applications to linear elliptic periodic homogenization.

(ii)
It is easy to verify that the linear operators $A_\ve$ do not converge for $\ve \to 0$ in the uniform operator norm in ${\cal L}(W_0^{1,p}(\Omega;\R^n);W^{-1,p}(\Omega;\R^n))$ for certain $p\ge 2$, in general (see \cite[Remark 8.4]{Ci} and Lemma
\ref{Klemma} below). This is the reason why the classical implicit function theorem is not directly applicable to the boundary value problem \reff{ODE}, in general.

But the sequences of matrix functions 
$[a_{ij}(\cdot/\ve_1)],[a_{ij}(\cdot/\ve_2)],\ldots$ H-converge (in the sense F. Murat and L. Tartar) to the matrix $[\hat a_{ij}]$ for any choice of a tending to zero sequence $\ve_1,\ve_2,\ldots>0$ 
(cf. Theorem \ref{Shentheorem} below), and this allows to apply a generalized implicit function theorem, namly Theorem \ref{ift}.
\end{remark}
Finally, we introduce a 
nonlinear operator $B:
C(\overline\Omega;\R^n)\to W^{-1,2}(\Omega;\R^n)$ by
\bee
\label{Fdef}
\langle B(u),\vp\rangle_{1,2}:=
\int_\Omega\Big(b_i^\al(x,u(x))\partial_{x_i}\vp^\al(x)+b^\al(x,u(x))\vp^\al(x)\Big)dx
\mbox{ for all } \vp \in W_0^{1,2}(\Omega;\R^n).
\ee
Here and in what follows we consider the function space $C(\overline\Omega;\R^n)$ with the norm $\|\cdot\|_\infty$ (defined in \reff{infdef}), as usual.
Remark that
for any $p>2$
the operator $B$ can be considered as a map from $W^{1,p}(\Omega;\R^n)$ into $W^{-1,p}(\Omega;\R^n)$
because the space $W^{1,p}(\Omega;\R^n)$
is continuously embedded into the space
$C(\overline\Omega;\R^n)$ (because the dimension of $\Omega$ is two).
Because of assumption \reff{diffass} we have that  
the nonlinear operator $B$ is $C^1$-smooth from $C(\overline\Omega;\R^n)$ into $W^{-1,2}(\Omega;\R^n)$,
and
$$
\langle B'(u)v,\vp\rangle_{1,2}:=
\int_\Omega\Big(\partial_{u^\gamma}b_i^\al(x,u(x))\partial_{x_i}\vp^\al(x)+\partial_{u^\gamma}b^\al(x,u(x))\vp^\al(x)\Big)v^{\gamma}(x)dx,
$$
and for all $u \in C(\overline\Omega;\R^n)$ and all $p \ge 2$ we have
\bee
\label{Bprop}
\lim_{\|v\|_\infty\to 0}
\|(B'(u+v)-B'(u))w\|_{-1,p}=0 \mbox{ uniformly with respect to } \|w\|_{\infty} \le 1.
\ee
Moreover,
using the notation \reff{Avedef} and \reff{Fdef} we get
\bee
\label{unulldef}
A_0u_0+B(u_0)=0.
\ee
In particular, from Theorem \ref{Gro} follows that
\bee 
\label{usmooth}
u_0 \in W_0^{1,p_1}(\Omega;\R^n).
\ee
Further, it is well-known (cf., e.g. \cite[Chapter 2.2]{Shen}) that there exists $p_2>2$ such that
\bee 
\label{vsmooth}
v^{\al \beta}_i\in W_{\rm loc}^{1,p_2}(\R^2).
\ee

Now we introduce the abstract setting of  Corollary \ref{cor} for the boundary value problem \reff{ODE}.  We take 
$p_0$ from the assumption in Theorem \ref{main} (ii),
$p_1$ from  Theorem \ref{Gro} and \reff{usmooth} and $p_2$ from \reff{vsmooth} and fix $p$ and $q$ as follows:
$$
2<p\le \min\{p_0,p_1,p_2\},
\; q:=\frac{p}{p-1}.
$$
The Banach spaces $U$ and $V$  and their norms are defined by
\begin{eqnarray*}
&&U:=W_0^{1,p}(\Omega;\R^n),\;
V:=W^{-1,p}(\Omega;\R^n)=W_0^{1,q}(\Omega;\R^n)^*,\\
&&\|\cdot\|_U:=\|\cdot\|_{1,p},\;
\|\cdot\|:=\|\cdot\|_\infty,\;
\|\cdot\|_V:=\|\cdot\|_{-1,p}.
\end{eqnarray*}
Because the space dimension is two, the assumption \reff{weaker} of Corollary \ref{cor} is satisfied in this setting.
Further, the $C^1$-smooth operators $F_\ve:U \to V$ of Theorem \ref{ift} are defined by 
$$
F_\ve(u):=A_\ve u+B(u).
$$
With these choices a  vector function $u$ is a weak solution to the boundary value problem \reff{ODE} if and only if $u$ belongs to the function space $U$ and satisfies the operator equation
$F_\ve(u)=0$.
Here we used Theorem \ref{Gro} again.
And finally, we define the exceptional approximate solution $u_0\in U$ of Corollary \ref{cor} to be the solution $u_0$ of the homogenized boundary value problem \reff{hombvp}, which is given in Theorem \ref{main}.

In order to prove Theorem \ref{main} we have to verify the conditions \reff{coerza} and \reff{Fas1}
in the setting introduced above, i.e. that there exist $\ve_0>0$ and $c>0$ such that
\bee
\label{coerz1}
\left.
\begin{array}{l}
\mbox{for all $\ve \in (0,\ve_0]$ the operators $A_\ve+B'(u_0)$
are Fredholm}\\
\mbox{of index zero from }
W_0^{1,p}(\Omega;\R^n) \mbox{ into }
W^{-1,p}(\Omega;\R^n),
\mbox{ and }\\ \|(A_\ve+B'(u_0))u\|_{-1,p} \ge c \|u\|_{1,p} \mbox{ for all } u \in W_0^{1,p}(\Omega;\R^n),
\end{array}
\right\}
\ee
and that
\bee
\label{Fas2}
\lim_{\|u\|_\infty \to 0}
\|(B'(u_0+u)-B'(u_0))v\|_{-1,p}=0
 \mbox{ uniformly with respect to } \|v\|_{1,p} \le 1,
\ee
and we have to construct a family $\bar u_\ve$ such that the properties \reff{lim} and \reff{newconv} are
satisfied in the setting introduced above, i.e. 
\bee
\label{newconv1}
\|\uu-u_0\|_\infty \to 0 \mbox{ for } 
\ve \to 0,
\ee
and
\bee
\label{lim1}
\|A_\ve\uu+B(\uu)\|_{-1,p}=o(1) \mbox{ for } \ve \to 0,
\ee
and, for proving assertion (ii) of Theorem \ref{main}, we have to verify that
\bee
\label{lim2}
\|\uu-u_0\|_\infty+\|A_\ve\uu+B(\uu)\|_{-1,p}=O(\ve^{1/p}) \mbox{ for } \ve \to 0.
\ee

Condition \reff{Fas2} is true because of \reff{Bprop}. Moreover, if \reff{newconv1} is proved, then
\reff{unulldef} and the continuity of $B$ imply that $\|B(\uu)-B(u_0)\|_{-1,p}=\|B(\uu)+A_0u_0\|_{-1,p}
\to 0$ for $\ve \to 0$, and then
the conditions \reff{lim1} and \reff{lim2}
are equivalent to
\bee
\label{lim5}
\|A_\ve\uu-A_0u_0\|_{-1,p}=o(1) \mbox{ for } \ve \to 0,
\ee
and
\bee
\label{lim6}
\|\uu-u_0\|_\infty+\|A_\ve\uu-A_0u_0\|_{-1,p}=O(\ve^{1/p}) \mbox{ for } \ve \to 0.
\ee
Hence, we have to verify \reff{coerz1}, and for proving assertion (i) of Theorem \ref{main} 
we have to
construct a family $\uu \in W_0^{1,p}(\Omega;\R^n)$ with  \reff{newconv1} and \reff{lim5}, and for proving assertion (ii) of Theorem \ref{main} we have construct a family $\uu \in W_0^{1,p}(\Omega;\R^n)$ with 
\reff{lim6}. This is what we are going to do below.

\subsection{Construction of approximate solutions with \reff{newconv1} and \reff{lim5}}
\label{sub1}
In this subsection we will construct a family $\uu\in  W_0^{1,p}(\Omega;\R^n)$ 
with 
\reff{newconv1} and \reff{lim5}.
For that we will do some calulations which are well-known in homogenization theory (cf., e.g. \cite[Chapter 3.2]{Shen}), and some estimates which seem to be new. 

For $\ve>0$ we set
$$
\Omega_\ve:=\left\{x\in \Omega:\; 
\inf_{y \in \partial\Omega}
\left((x_1-y_1)^2+(x_2-y_2)^2\right)<\ve^2\right\}.
$$
It follows that
\bee
\label{Omest}
|\Omega_\ve|=O(\ve)
\mbox{ for } \ve \to 0,
\ee
where $|\Omega_\ve|$ is the two-dimensional Lebesque measure
of $\Omega_\ve$.
Further, we take a family $\eta_\ve$ of
cut-of functions of size $\ve$, i.e. of $C^\infty$-functions $\Omega\to \R$ such that
\bee
\label{etaest}
\left.
\begin{array}{l}
0 \le \eta_\ve(x)\le 1 \mbox{ for } x \in \Omega,\\
\eta_\ve(x)=1\mbox{ for } x \in  \Omega\setminus \Omega_{2\ve},\\
\eta_\ve(x)=0\mbox{ for } x \in \Omega_{\ve},\\
\sup\left\{\ve\,|\partial_{x_i}\eta_\ve(x)|:\; \ve>0, \;x \in \Omega,\; i=1,2\right\}<\infty.
\end{array}
\right\}
\ee
Finally, we take a mollifier function, i.e. a $C^\infty$-function $\rho:\R^2\to \R$ such that
$$
\rho(x)\ge 0 \mbox{ and } \rho(x)=\rho(-x) \mbox{ for all } x \in \R^2,\;
\rho(x)=0 \mbox{ for } x_1^2+x_2^2 \ge 1,
\int_{\R^2}\rho(x)dx=1,
$$
and for $\delta>0$ we define linear smoothing operators $S_\delta: L^1(\Omega) \to C^\infty(\R^2)$ by
$$
[S_\delta u](x):=\int_{\Omega}\rho_\delta(x-y) u(y)dy
\mbox{ with }
\rho_\delta(x):=\rho(x/\delta)/\delta^2.
$$

\begin{lemma}
\label{Slemma}
(i) For all $r \ge 1$ and $u \in L^r(\Omega)$ we have
\bee
\label{Sconv}
\lim_{\delta\to 0}\int_\Omega|u(x)-[S_\delta u](x)|^rdx=0
\ee

(ii) For all $r \ge 1$ there exists $c_r>0$ such that for all $\delta>0$
and $u \in L^r(\Omega)$ we have
\begin{eqnarray}
\label{pest1}
\int_\Omega\left|[\partial_{x_i}S_\delta u](x)\right|^rdx &\le& \frac{c_r}{\delta^r}
\int_\Omega\left|u(x)\right|^rdx
\mbox{ for } i=1,2,\\
\label{pest3}
\sup_{x \in \Omega}\left|[S_\delta u](x)\right|^r&\le& 
 \frac{c_r}{\delta^{2}}
\int_\Omega\left|u(x)\right|^rdx.
\end{eqnarray}
\end{lemma}
{\bf Proof }
Assertion (i) is proved e.g. in
\cite[Lemma 1.1.1]{Barbu}.

In order to prove assertion (ii) take $\delta>0$, $r,s>1$ with $1/r+1/s=1$, and take $u \in L^r(\Omega)$.
Then the H\"older inequality implies that for all $x\in \Omega$ we have 
$$
|[S_\delta u](x)|=\left|\int_{\Omega} u(y)\rho_\delta(x-y)^{1/r}\rho_\delta(x-y)^{1/s}dy\right|\le\left(\int_{\Omega}|u(y)|^r\rho_\delta(x-y)dy\right)^{1/r}.
$$
Here we used that
$
\int_{\R^2}\rho_\delta(x-y)dy
=
\int_{\R^2}\rho(z)dz=1.
$
It follows that
$$
|[S_\delta u](x)|^r
\le \frac{1}{\delta^2}
\int_{\Omega}|u(y)|^r\rho((x-y)/\delta)dy
\le \mbox{const } \frac{1}{\delta^2}\int_{\Omega}|u(y)|^rdy,
$$
where the constant does not depend on $\delta$ and $u$.
Hence, \reff{pest3} is proved.

Further, because of
$
\int_{\R^2}|\partial_{x_i}\rho_\delta(x-y)|dx=\frac{1}{\delta^3}\int_{\R^2}|\partial_{x_i}\rho((x-y)/\delta)|dy=
\frac{1}{\delta}\int_{\R^2}|\partial_{x_i}\rho(z)|dz
$
we have that
\begin{eqnarray*}
|[\partial_{x_i}S_\delta u](x)|&\le&
\left(\int_{\Omega}| u(y)|^r|\partial_{x_i}\rho_\delta(x-y)|dy\right)^{1/r}
\left(\int_{\Omega}|\partial_{x_i}\rho_\delta(x-y)|dy\right)^{1/s}\\
&\le& \mbox{ const }
\frac{1}{\delta^{1/s}}\left(\int_{\Omega}|u(y)|^r|\partial_{x_i}\rho_\delta(x-y)|dy\right)^{1/r}
\end{eqnarray*}
and, hence,
\begin{eqnarray*}
\int_\Omega\left|[\partial_{x_i}S_\delta u](x)\right|^rdx&\le& \mbox{ const }
\frac{1}{\delta^{r/s}}\int_{\Omega}| u(y)|^r\int_{\Omega}|\partial_{x_i}\rho_\delta(x-y)|dx\;dy\\
&\le&\mbox{ const }
\frac{1}{\delta^r}\int_{\Omega}|u(y)|^rdy,
\end{eqnarray*}
where, again, the constants do not depend on $\delta$ and $u$.
Hence, \reff{pest1} is proved also.
\qed\\

Because of
\reff{usmooth} and \reff{vsmooth} we can define
the needed family $\uu \in W_0^{1,p}(\Omega;\R^n)$, which should satisfy
\reff{newconv1} and \reff{lim5},
as follows:
\bee
\label{barudef1}
\uu^\al(x):=
u^\al_0(x)+\ve \eta_\ve(x) [S_{\delta_\ve}\partial_{x_k}u_0^\gamma](x) v_k^{\al \gamma}(\x)
\mbox{ with } \delta_\ve:=\ve^{1/4}.
\ee
\begin{remark}
\label{choice}
In fact we need a choice of $\delta_\ve$ in \reff{barudef1} such that $\delta_\ve$ tends to zero for $\ve \to 0$, but this convergence should be sufficiently slow. More exactly, we need that
$$
\lim_{\ve \to 0}\left(\|S_{\delta_\ve}\partial_{x_k}u_0-\partial_{x_k}u_0\|_{L^p}+
\ve \|S_{\delta_\ve}\partial_{x_k}u_0\|_\infty+\ve
\|\partial_{x_j}S_{\delta_\ve}
\partial_{x_k}u_0\|_{L^p}\right)=0
$$
(see the estimates \reff{Kest} and 
\reff{es1}-\reff{Sdeltaest} below).
For example, the choices $\delta_\ve=\ve^r$ with $0<r<1/2$ or
$\delta_\ve=-1/\ln \ve$ are appropriate.
\end{remark}

From \reff{pest3} and from definition 
\reff{barudef1} it follows that
\bee
\label{Kest}
\|\uu-u_0\|_\infty \le \mbox{const }\frac{\ve}{\delta_\ve^{2/p}}\|u_0\|_{1,p}\le \mbox{const }\ve^{(2p-1)/2p},
\ee
where the constants do not depend on $\ve$.
Hence,
\reff{newconv1} is satisfied.
In order to prove \reff{lim5} we use the following lemma:
\begin{lemma}
\label{Klemma}
We have
$\|A_\ve\uu-A_0u_0\|_{-1,p}\to 0$
for $\ve \to 0$.
\end{lemma}
{\bf Proof }
For $\al,\beta=1,\ldots,n$ and $i,j,k=1,2$
we define
$\Z^2$-periodic functions $b_{ij}^{\al \beta}\in L_{\rm loc}^{p}(\R^2)$ and $c_{ij}^{\al \beta}\in W_{\rm loc}^{2,p}(\R^2)$
and 
$\phi_{ijk}^{\al \beta}\in W_{\rm loc}^{1,p}(\R^2)$ 
(the functions $\phi_{ijk}^{\al \beta}$  sometimes are called dual or flux correctors)  
by
\bee
\label{bdef}
b_{ij}^{\al \beta}(y):=a_{ij}^{\al \beta}(y)
+a^{\al \gamma}_{ik}(y)
\partial_{y_k}v_j^{\gamma \beta}(y)-\hat a^{\al \beta}_{ij}
\ee
and
\bee
\label{cdef}
\Delta c_{ij}^{\al \beta}(y)=b_{ij}^{\al \beta}(y),\; \int_{[0,1]^2}c_{ij}^{\al \beta}(y)dy=0
\ee
and
\bee
\label{phidef}
\phi_{ijk}^{\al \beta}(y):=
\partial_{y_i}c_{jk}^{\al \beta}(y)
-\partial_{y_j}c_{ik}^{\al \beta}(y).
\ee
From \reff{Anullcoef} and \reff{bdef} follows that $\int_{[0,1]^2}b^{\al \beta}_{ij}(y)dy=0$,
therefore problem \reff{cdef} is uniquely strongly solvable with respect to $c^{\al \beta}_{ij}$. Further,
from \reff{cell} follows that $\partial_{y_i}b_{ij}^{\al \beta}=0$.
Hence, \reff{cdef} implies that $\partial_{y_i}c_{ij}^{\al \beta}=0$.
Therefore \reff{cdef} and \reff{phidef} yield that
\bee
\label{phiprop}
\partial_{y_i}\phi_{ijk}^{\al \beta}=b_{jk}^{\al \beta} 
\mbox{ and }
\phi_{ijk}^{\al \beta}=-\phi_{kji}^{\al \beta}.
\ee
Using \reff{phiprop} we get
\bee
\label{ibp}
\ve \partial_{x_k}\left(\phi_{ijk}^{\al \beta}(\x)\partial_{x_i}\vp^\beta(x)\right)=
b_{ij}^{\al \beta}(\x)\partial_{x_i}\vp^\beta(x)
\mbox{ for all } \vp \in C^\infty(\Omega;\R^n)
\ee 
(this is \cite[formula (3.1.5)]{Shen}).

Now, we 
insert the definition \reff{barudef1} of $\uu$ into $\langle A_\ve\uu-A_0u_0,\vp\rangle_{1,q}$
with $\vp \in C^\infty(\Omega;\R^n)$
and calculate as follows:
\begin{eqnarray}
&&
\langle A_\ve\uu-A_0u_0,\vp\rangle_{1,q}
\nonumber\\
&&=\int_\Omega\left(
a^{\al \beta}_{ij}(\x)\partial_{x_j}\left(u_0^\beta+\ve\eta_\ve
[S_{\delta_\ve}\partial_{x_k}u_0^\gamma]
v_k^{\beta \gamma}(\x)
\right)                                    
-\hat a^{\al \beta}_{ij}\partial_{x_j}
u_0^\beta\right)
\partial_{x_i}\vp^\al dx\nonumber\\
&&=
\int_\Omega\left(
\left(a^{\al \beta}_{ij}(\x)
-\hat a^{\al \beta}_{ij}\right)
\partial_{x_j}u_0^\beta
+a^{\al \beta}_{ij}(\x)
\eta_\ve
[S_{\delta_\ve}\partial_{x_k}u_0^\gamma]
\partial_{y_j}v_k^{\beta \gamma}(\x)
\right)                                    
\partial_{x_i}\vp^\al dx\nonumber\\
&&\;\;\;\;\;\;+\ve\int_\Omega
a^{\al \beta}_{ij}(\x)\partial_{x_j}(
\eta_\ve(x)
[S_{\delta_\ve}\partial_{x_k}u_0^\gamma](x))
v_k^{\beta \gamma}(\x)
\partial_{x_i}\vp^\al(x)
dx\nonumber\\
&&=
\int_\Omega\left(
a^{\al \beta}_{ij}(\x)
-\hat a^{\al \beta}_{ij}
+a^{\al \gamma}_{ik}(\x)
\partial_{y_k}v_j^{\gamma\beta}(\x)\right)\eta_\ve(x)
[S_{\delta_\ve}\partial_{x_j}u_0^\beta](x)
\partial_{x_i}\vp^\al(x)dx\nonumber\\
&&\;\;\;\;\;\;+\int_\Omega
\left(a^{\al \beta}_{ij}(\x)
-\hat a^{\al \beta}_{ij}\right)\left(
\partial_{x_j}u_0^\beta(x)-
\eta_\ve(x)
[S_{\delta_\ve}\partial_{x_j}u_0^\beta](x)\right)
\partial_{x_i}\vp^\al(x)dx\nonumber\\
&&\;\;\;\;\;\;+\ve\int_\Omega
a^{\al \beta}_{ij}(\x)\partial_{x_j}(
\eta_\ve(x)
[S_{\delta_\ve}\partial_{x_k}u_0^\gamma](x))
v_k^{\beta \gamma}(\x)
\partial_{x_i}\vp^\al(x)
dx.
\label{int}
\end{eqnarray}
We insert \reff{bdef} and \reff{ibp} into \reff{int}, integrate by parts and use that $\phi_{kij}^{\al \beta}(\x)\partial_{x_k}\partial_{x_i}\vp^\al(x)=0$
(cf. \reff{phiprop})
and get
\begin{eqnarray}
&&
\langle A_\ve\uu-A_0u_0,\vp\rangle_{1,q}
\nonumber\\
&&=\int_\Omega
b_{ij}^{\al \beta}(\x)
\eta_\ve(x)
[S_{\delta_\ve}\partial_{x_j}u_0^\beta](x)
\partial_{x_i}\vp^\al(x)
dx\nonumber\\
&&\;\;\;\;\;\;+\int_\Omega
\left(a^{\al \beta}_{ij}(\x)
-\hat a^{\al \beta}_{ij}\right)\left(
\partial_{x_j}u_0^\beta(x)-
\eta_\ve(x)
[S_{\delta_\ve}\partial_{x_j}u_0^\beta](x)\right)
\partial_{x_i}\vp^\al(x)dx\nonumber\\
&&\;\;\;\;\;\;+\ve\int_\Omega
a^{\al \beta}_{ij}(\x)\partial_{x_j}(
\eta_\ve(x)
[S_{\delta_\ve}\partial_{x_k}u_0^\gamma](x))
v_k^{\beta \gamma}(\x)
\partial_{x_i}\vp^\al(x)
dx\nonumber\\
&&=\ve\int_\Omega\left(-\phi_{ijk}^{\al \beta}(\x)+
a^{\al\gamma}_{ik}(\x) v_j^{\gamma\beta}(\x)
\right)
\partial_{x_k}(\eta_\ve(x)
[S_{\delta_\ve}\partial_{x_j}u_0^\beta](x))
\partial_{x_i}\vp^\al(x)dx\nonumber\\
&&\;\;\;\;\;\;+\int_\Omega
\left(a^{\al \beta}_{ij}(\x)
-\hat a^{\al \beta}_{ij}\right)\left(
\partial_{x_j}u_0^\beta(x)-
\eta_\ve(x)
[S_{\delta_\ve}\partial_{x_j}u_0^\beta](x)\right)
\partial_{x_i}\vp^\al(x)dx.
\label{threeint}
\end{eqnarray}
Remark that no boundary integrals appeared after the integration by parts because
of the cut-off functions $\eta_\ve$ (no matter if the test function $\vp$ vanishes on $\partial \Omega$ or not).

Let us estimate the the right-hand side of \reff{threeint}. Because of  \reff{usmooth} and \reff{vsmooth} and the H\"older inequality
we have
\begin{eqnarray}
&&\left|\ve\int_\Omega\left(-\phi_{kij}^{\al \beta}(\x)+
a^{\al \gamma}_{ik}(\x) v_j^{\gamma \beta}(\x)
\right)\partial_{x_k}\eta_\ve(x)
[S_{\delta_\ve}\partial_{x_j}u_0^\beta](x)
\partial_{x_i}\vp^\al(x)dx\right|\nonumber\\
&&\le\mbox{const }\ve\left(\sum_{i=1}^2
\int_{\Omega_\ve}
|\partial_{x_i}\eta_\ve(x)|^pdx\right)^{1/p}
\sum_{\beta=1}^n\sum_{j=1}^2\|S_{\delta_\ve}\partial_{x_j}
u_0^\beta\|_\infty
\|\vp\|_{1,q}\nonumber\\
&&\le\mbox{const } \frac{\ve^{1/p}} 
{\delta_\ve^{2/p}}\|u_0\|_{1,p}
\|\vp\|_{1,q}
\le
\mbox{const }\ve^{1/2p}
\|\vp\|_{1,q}
\label{es1}
\end{eqnarray}
(here we used \reff{Omest}, \reff{etaest} and \reff{pest3})
and
\begin{eqnarray}
&&\left|\ve\int_\Omega\left(-\phi_{kij}^{\al \beta}(\x)+
a^{\al \gamma}_{ik}(\x) v_j^{\gamma \beta}(\x)
\right)\eta_\ve(x)
[\partial_{x_k}
S_{\delta_\ve}\partial_{x_j}u_0^\beta](x)
\partial_{x_i}\vp^\al(x)dx\right|\nonumber\\
&&\le\mbox{const }
\frac{\ve}{\delta_\ve}
\|u_0\|_{1,p}
\|\vp\|_{1,q}\le\mbox{const }
\ve^{3/4}
\|\vp\|_{1,q}
\label{es2}
\end{eqnarray}
(here we used \reff{pest1})
and
\begin{eqnarray}
&&\left|
\int_\Omega
\left(a^{\al \beta}_{ij}(\x)
-\hat a^{\al \beta}_{ij}\right)\left(
1-
\eta_\ve(x)\right)
[S_{\delta_\ve}\partial_{x_j}u_0^\beta](x)
\partial_{x_i}\vp^\al(x)dx\right|\nonumber\\
&&\le\mbox{const}\left(
\int_{\Omega_\ve}
|1-\eta_\ve(x)|^pdx\right)^{1/p}
\sum_{\beta=1}^n\sum_{j=1}^2\|[S_{\delta_\ve}\partial_{x_j}
u_0^\beta]\|_\infty
\|\vp\|_{1,q}\nonumber\\
&&\le\mbox{const } \frac{\ve^{1/p}}
{\delta_\ve^{2/p}}\|u_0\|_{1,p}
\|\vp\|_{1,q}
\le\mbox{const } \ve^{1/2p}
\|\vp\|_{1,q}
\label{es3}
\end{eqnarray}
(here we used \reff{Omest} and \reff{pest3}),
where the constants do not depend on $\ve$
and $\vp$.
Further, we have
\begin{eqnarray}
\label{Sdeltaest}
&&\left|
\int_\Omega
\left(a^{\al \beta}_{ij}(\x)
-\hat a^{\al \beta}_{ij}\right)\left(
\partial_{x_j}u_0^\beta(x)-
[S_{\delta_\ve}\partial_{x_j}u_0^\beta](x)\right)
\partial_{x_i}\vp^\al(x)dx\right|\nonumber\\
&&\le \mbox{const}\left(\sum_{\beta=1}^n\sum_{j=1}^2\int_\Omega\left|\partial_{x_j}u_0^\beta(x)-
[S_{\delta_\ve}\partial_{x_j}u_0^\beta](x)\right|^{p}dx\right)^{1/p}
\|\vp\|_{1,q},
\end{eqnarray}
where the constant does not depend on $\ve$
and $\vp$ again.
But the right-hand side of \reff{Sdeltaest} is $o(1)$ for $\ve \to 0$ uniformly with respect to $\|\vp\|_{1,q}\le 1$ (cf. \reff{Sconv}).
Hence, the lemma is proved.
\qed
\subsection{Construction of approximate solutions with 
\reff{lim6}}
\label{sub2}
In this subsection we will  construct a family $\uu\in  W_0^{1,p}(\Omega;\R^n)$ 
with 
\reff{lim6} under the assumption that
$u_0 \in W_0^{2,p_0}(\Omega;\R^n)$ with certain $p_0>2$.
Therefore
in the definition of the
family $\uu$  
we do not need the smoothing operators $S_\delta$, i.e. this time we define
\bee
\label{barudef2}
\uu^\al(x):=u^\al_0(x)+\ve \eta_\ve(x) \partial_{x_k}u_0^\gamma(x) v_k^{\al \gamma}(\x),
\ee
and because of \reff{vsmooth} 
we have that
$\|\uu-u_0\|_\infty=O(\ve)$ for $\ve \to 0$.

In order to verify \reff{lim6} 
it remains to show that $\|A_\ve\uu-A_0u_0\|_{-1,p}=O(\ve^{1/p})$ for $\ve \to 0$.
To show this 
we proceed  as in \reff{threeint}:
\begin{eqnarray*}
&&\int_\Omega\Big(
a^{\al \beta}_{ij}(\x)\partial_{x_j}
\uu^\beta(x)
-\hat a^{\al \beta}_{ij}\partial_{x_j}
u^\beta_0(x)\Big)
\partial_{x_i}\vp^\al(x)dx\nonumber\\
&&=
\int_\Omega\left(
a^{\al \beta}_{ij}(\x)\partial_{x_j}\left(u^\beta_0+\ve\eta_\ve
\partial_{x_k}u_0^\gamma
v_k^{\beta \gamma}(\x)
\right)                                    
-\hat a^{\al \beta}_{ij}\partial_{x_j}
u^\beta_0\right)
\partial_{x_i}\vp^\al dx\nonumber\\
&&=\ve\int_\Omega\left(-\phi_{ijk}^{\al \beta}(\x)+
a^{\al\gamma}_{ik}(\x) v_j^{\gamma\beta}(\x)
\right)
\partial_{x_k}(\eta_\ve(x)
\partial_{x_j}u_0^\beta(x))
\partial_{x_i}\vp^\al(x)dx\nonumber\\
&&\;\;\;\;\;\;+\int_\Omega
\left(a^{\al \beta}_{ij}(\x)
-\hat a^{\al \beta}_{ij}\right)\left(
\partial_{x_j}u_0^\beta(x)-
\eta_\ve(x)
\partial_{x_j}u_0^\beta(x)\right)
\partial_{x_i}\vp^\al(x)dx.
\end{eqnarray*}
And as in \reff{es1}-\reff{es3} one estimates as follows:
\begin{eqnarray*}
&&\left|\ve\int_\Omega\left(-\phi_{ijk}^{\al \beta}(\x)+
a^{\al \gamma}_{ik}(\x) v_j^{\gamma \beta}(\x)
\right)\partial_{x_k}\eta_\ve(x)
\partial_{x_j}u_0^\beta(x)
\partial_{x_i}\vp^\al(x)dx\right|\nonumber\\
&&\le\mbox{const }\ve\left(\sum_{i=1}^2
\int_{\Omega_\ve}
|\partial_{x_i}\eta_\ve(x)|^pdx\right)^{1/p}
\sum_{\beta=1}^n\sum_{j=1}^2\|\partial_{x_j}
u^\beta_0\|_\infty
\|\vp\|_{1,q}
\le
\mbox{const }\ve^{1/p}\|\vp\|_{1,q}
\end{eqnarray*}
and
$$
\left|\ve\int_\Omega\left(-\phi_{ijk}^{\al \beta}(\x)+
a^{\al \gamma}_{ik}(\x) v_j^{\gamma \beta}(\x)
\right)\eta_\ve(x)
\partial_{x_k}
\partial_{x_j}u_0^\beta(x)
\partial_{x_i}\vp^\al(x)dx\right|\nonumber\le\mbox{const }
\ve
\|\vp\|_{1,q}
$$
and
\begin{eqnarray*}
&&\left|
\int_\Omega
\left(a^{\al \beta}_{ij}(\x)
-\hat a^{\al \beta}_{ij}\right)\left(
1-
\eta_\ve(x)\right)
\partial_{x_j}u_0^\beta(x)
\partial_{x_i}\vp^\al(x)dx\right|\nonumber\\
&&\le\mbox{const}\left(
\int_{\Omega_\ve}
|1-\eta_\ve(x)|^pdx\right)^{1/p}
\|\vp\|_{1,q}
\le\mbox{const } \ve^{1/p}
\|\vp\|_{1,q},
\end{eqnarray*}
where the constants do not depend on $\ve$ and $\vp$.
Hence, \reff{lim6} is proved.

\subsection{Verification of  \reff{coerz1}}
\label{sub3}
The linear operators $A_\ve$ are isomorphisms from $W_0^{1,p}(\Omega;\R^n)$ onto $W^{-1,p}(\Omega;\R^n)$ (cf. Theorem~\ref{Gro}), and the linear operators $B'(u_0)$ are bounded from $C(\overline\Omega;\R^n)$ into
$W^{-1,p}(\Omega;\R^n)$
and, hence, compact from $W_0^{1,p}(\Omega;\R^n)$ into $W^{-1,p}(\Omega;\R^n)$. Hence, condition \reff{coerz1} is satisfied if there exists $\ve_0>0$
such that
$$
\inf\left\{\|(A_\ve+B'(u_0))u\|_{-1,p}:\;
\ve \in (0,\ve_0], u \in W_0^{1,p}(\Omega;\R^n), \|u\|_{1,p}= 1\right\}>0.
$$

Suppose that this is not true. Then there exist sequences $\ve_1,\ve_2,\ldots>0$ and $u_1,u_2,\ldots \in W_0^{1,p}(\Omega;\R^n)$ such that $\ve_l\to 0$ for $l \to \infty$ and
\bee
\label{conve1}
\lim_{l\to \infty}\|(A_{\ve_l}+B'(u_0))u_l\|_{-1,p}=0,
\ee
but
\bee
\label{normone1}
\|u_l\|_{1,p}= 1 \mbox{ for all } l.
\ee
Because $W_0^{1,p}(\Omega;\R^n)$ is reflexive and because it is compactly embedded into  $C(\overline\Omega;\R^n)$, without loss of generality we may assume that there exists $u_*\in W_0^{1,p}(\Omega;\R^n)$ such that 
\bee
\label{infconv}
u_l\rightharpoonup u_* \mbox{ for $l\to \infty$
weakly in }
W^{1,p}(\Omega;\R^n)
\mbox{ and }
\lim_{l \to \infty}\|u_l-u_*\|_\infty=0.
\ee
From \reff{infconv} follows that 
$\|B'(u_0)(u_l-u_*)\|_{-1,p}\to 0$ for
$l\to \infty$. Therefore \reff{conve1} implies that
\bee
\label{conve2}
\lim_{l\to \infty}\|A_{\ve_l}u_l+B'(u_0)u_*\|_{-1,p}=0.
\ee

Now we use the following theorem, which is well-known in periodic homogenization theory for linear elliptic equations and  systems with $L^\infty$-coefficients (see, e.g. \cite[Lemma 8.6]{Che}, \cite[Theorem 2.3.2]{Shen}). Its proof is based on the Div-Curl Lemma.
\begin{theorem}
\label{Shentheorem}
Let be given $\tilde u_* \in W^{1,2}(\Omega;\R^n)$ and $f_*\in W^{-1,2}(\Omega;\R^n)$ and sequences $\ve_1,\ve_2,\ldots>0$ and $\tilde u_1,\tilde u_2,\ldots \in 
W^{1,2}(\Omega;\R^n)$ such that
$$
\tilde u_l\rightharpoonup \tilde u_* \mbox{ for $l\to \infty$ weakly in 
$W^{1,2}(\Omega;\R^n)$ and }
\lim_{l \to \infty}\left(\ve_l+\|A_{\ve_l}\tilde u_l-f_*\|_{-1,2}\right)=0.
$$
Then $A_0\tilde u_*=f_*$.
\end{theorem}

Because of \reff{infconv}, \reff{conve2} and Theorem \ref{Shentheorem} it follows that
$
(A_0+B'(u_0))u_*=0,
$
i.e. that
$u_*$ is a weak solution
to the linearized boundary value problem
\reff{linbvp}. Hence, by assumption of Theorem \ref{main}, we get that $u_*=0$.
Therefore
\reff{conve2} implies that $\|A_{\ve_l}u_l\|_{-1,p}\to 0$ for $l \to \infty$.
But this contradicts to \reff{Ainvert} and \reff{normone1}.

\section{Nonlinear mixed boundary conditions}
\label{sec4}
\setcounter{equation}{0}
\setcounter{theorem}{0}

In this section we generalize the results of Theorem \ref{main} to the case of boundary value problems of the type
\bee
\label{Robin}
\left.
\begin{array}{l}
\partial_{x_i}\Big(a^{\al \beta}_{ij}(\x)
\partial_{x_j}u^\beta(x)
+b_i^\al(x,u(x))
\Big)=b^\al(x,u(x)) \mbox{ for } x \in \Omega,\\
\Big(a^{\al \beta}_{ij}(\x)\partial_{x_j}u^\beta(x)
+b_i^\al(x,u(x))\Big)\nu_i(x)=b_0^\al(x,u(x))
\mbox{ for } x \in \Gamma,\\
u^\al(x)=0 \mbox{ for } x \in \partial\Omega\setminus\Gamma.
\end{array}
\right\}
\al=1,\ldots,n.
\ee
Here $\nu=(\nu_1,\nu_2):\partial \Omega\to \R^2$ is the outer unit normal vector field on the boundary $\partial \Omega$,
and $u \in \R^n \mapsto b_0(\cdot,u)\in L^\infty(\Gamma;\R^n)$ is $C^1$-smooth.

The reasons, why similar to Theorem \ref{main} results are true also for the boundary value problem \reff{Robin}, are easily to explain:
First, Theorem \ref{Gro} is true (in an appropriate reformulation) also for mixed boundary conditions. Second, in 
\reff{threeint}
 we did not need that the test functions $\vp$ vanish on $\partial \Omega$. And third, in Theorem \ref{Shentheorem} no boundary conditions are imposed on the functions $\tilde u_l$.
 
 We assume that $\Omega$ is a bounded Lipschitz domain in $\R^n$ and that the Robin condition boundary part $\Gamma$ 
is a subset of $\partial \Omega$ such that
$$
\mbox{either $\Gamma=\emptyset$ 
or $\Gamma=\partial\Omega$ or the relative boundary of  
$\Gamma$ in $\partial \Omega$ consists of finitely many points.}
$$
Our assumptions concerning 
the diffusion coefficients $a^{\al \beta}_{ij}$ and the drift and reaction functiond $b_i^\al$ and $b^\al$ are as in Section \ref{sec:1} with one exception: If $\Gamma =\partial \Omega$, i.e. in the case of pure Robin boundary conditions, in place of condition \reff{alass} we assume the stronger condition \reff{Leg}.

The homogenized diffusion coefficients $\hat a_{ij}$ are as in 
Section~\ref{sec:1}, i.e. defined in 
\reff{Anullcoef}
via the cell problems \reff{cell},
and the homogenized boundary value problem, corresponding to \reff{Robin}, is
\bee
\label{homRobin}
\left.
\begin{array}{l}\partial_{x_i}\Big(\hat a^{\al \beta}_{ij}
\partial_{x_j}u^\beta(x)
+b_i^\al(x,u(x))\Big)
=b^\al(x,u(x)) \mbox{ for } x \in \Omega,\\
\Big(\hat a^{\al \beta}_{ij}\partial_{x_j}u^\beta(x)
+b_i^\al(x,u(x))\Big)
\nu_i(x)=b^\al_0(x,u(x))\mbox{ for } x \in \Gamma,\\
u^\al(x)=0 \mbox{ for } x \in \partial\Omega\setminus\Gamma.
\end{array}
\right\}
\al=1,\ldots,n.
\ee

A vector function $u \in W^{1,2}(\Omega;\R^n)
\cap C(\overline\Omega;\R^n)$ is called weak solution to the boundary value problem \reff{Robin} if it vanishes on $\partial \Omega \setminus \Gamma$ and if it satisfies the variational equation
\begin{eqnarray*}
&&\int_\Omega
\left(\Big(a^{\al \beta}_{ij}(\x)
\partial_{x_j}u^\beta(x)
+b_i^\al(x,u(x))\Big)
\partial_{x_i}\vp^\al(x)
+b^\al(x,u(x))\vp^\al(x)\right)dx\\
&&=\int_{\Gamma}
b_0^\al(x,u(x))\vp^\al(x)d\mu(x)
\mbox{ for all }
\vp \in W_\Gamma^{1,2}(\Omega;\R^n),
\end{eqnarray*}
where $\mu$  is the (one-dimensional) Lebesgue measure on $\partial\Omega$, and similarly for \reff{homRobin} and for its linearization
\bee
\label{linhomRobin}
\left.
\begin{array}{l}\partial_{x_i}\Big(
\hat a^{\al \beta}_{ij}
\partial_{x_j}u^\beta(x)
+\partial_{u^\gamma}b_i^\al(x,u_0(x))u^\gamma(x)
\Big)
=\partial_{u^\gamma}b^\al(x,u_0(x))u^\gamma(x) \mbox{ for } x \in \Omega,\\
\Big(\hat a^{\al \beta}_{ij}\partial_{x_j}u^\beta(x)
+\partial_{u^\gamma}b_i^\al(x,u_0(x))u^\gamma(x)
\Big)
\nu_i(x)=
\partial_{u^\gamma}
b_0^\al(x,u_0(x))u^\gamma(x)
\mbox{ for } x \in \Gamma,\\
u^\al(x)=0 \mbox{ for } x \in \partial\Omega\setminus\Gamma
\end{array}
\right\}
\ee
for $\al=1,\ldots,n$. Similar to Theorem \ref{main}, we get the following

\begin{theorem} 
\label{mainRob}
Let $u=u_0$ be a weak solution to \reff{homRobin} such that
\reff{linhomRobin}
does not have weak solutions $u\not=0$.
Then the following is true:

(i)
There exist $\ve_0>0$ and  $\delta>0$
such that for all $\ve \in (0,\ve_0]$ there exists exactly one 
weak solution $u=u_\ve$ 
to \reff{Robin} with $\|u-u_0\|_\infty \le \delta$. Moreover,
\bee
\label{estim1}
\|u_\ve-u_0\|_\infty\to 0 \mbox{
for }
\ve \to 0.
\ee

(ii) If 
$u_0 \in W^{2,p_0}(\Omega;\R^n)$
with certain $p_0>2$, then for all $p>2$ we have
\bee
\label{estim2}
\|u_\ve-u_0\|_\infty=O(\ve^{1/p})  \mbox{ for } \ve \to 0.
\ee
\end{theorem}

The proof of Theorem \ref{mainRob} is similar to that of Theorem \ref{main}. We indicate only the few small differences.
One has to apply 
Corollary \ref{cor}
again, but now in the following setting:
$$
U:=W_\Gamma^{1,p}(\Omega;\R^n),\;
V:=W_\Gamma^{1,q}(\Omega;\R^n)^*
\mbox{ with } 1/p+1/q=1,
$$ 
where $\|\cdot\|_U:=\|\cdot\|_{1,p}$ (the norm $\|\cdot\|_{1,p}$ is defined in \reff{1pnorm}, and
$W_\Gamma^{1,p}(\Omega;\R^n)$ is the closure with respect to this norm of the set of all $\vp \in C^\infty(\Omega;\R^n)$ with support in $\Omega \cup \Gamma$,
and $\|\cdot\|:=\|\cdot\|_\infty$,
and
$$
\|f\|_V:=\sup\{\langle f,\vp\rangle_{1,q,\Gamma}:\;\vp \in W_\Gamma^{1,q}(\Omega;\R^n),\; \|\vp\|_{1,q}\le 1\},
$$
where $\langle\cdot,\cdot\rangle_{1,q,\Gamma}:
W_\Gamma^{1,q}(\Omega;\R^n)^*\times
W_\Gamma^{1,q}(\Omega;\R^n)\to \R$ is the dual pairing, again.
The $C^1$-smooth operators $F_\ve:U \to V$ of Theorem \ref{ift} are defined by 
$F_\ve(u):=A_\ve u+B(u)$,
again, where 
the linear operators $A_\ve:W_\Gamma^{1,p}(\Omega;\R^n)\to W_\Gamma^{1,q}(\Omega;\R^n)^*$ and the nonlinearity $B:C(\overline\Omega;\R^n)\to W_\Gamma^{1,q}(\Omega;\R^n)^*$ 
are defined by
\begin{eqnarray*}
\langle A_\ve u,\vp\rangle_{1,q,\Gamma}&:=&
\int_\Omega \left(a_{ij}^{\al \beta}(\x)
\partial_{x_j}u^\beta(x)
\partial_{x_i}\vp^\al(x)+u^\al(x)\vp^\al(x)\right)dx,\\
\langle B(u),\vp\rangle_{1,q,\Gamma}&:=&
\int_\Omega\Big(b_i^\al(x,u(x))
\partial_{x_i}\vp(x)
+\left(b^\al(x,u(x))-
u^\al(x)\right)\vp^\al(x)\Big)dx\\
&&-
\int_{\Gamma} b_0^\al(x,u(x))\vp^\al(x)d\Gamma(x)
\end{eqnarray*}
for all $\vp \in W_\Gamma^{1,q}(\Omega;\R^n)$,
and similarly $A_0:W_\Gamma^{1,p}(\Omega;\R^n)\to W_\Gamma^{1,q}(\Omega;\R^n)^*$.

For proving the error estimates \reff{estim1} and \reff{estim2} we use again the families of approximate solutions \reff{barudef1} and \reff{barudef2}, respectively.

\begin{remark}
Let us draw attention to the following technical detail: 

In  Section \ref{sec3}, i.e. in the case of Dirichlet boundary conditions, the reason for using the cut-off functions $\eta_\ve$ in \reff{barudef1} and \reff{barudef2} is that the approximate solutions $\bar u_\ve$ should satisfy the Dirichlet boundary conditions. There 
the cut-off functions $\eta_\ve$ are not needed
to avoid  
boundary integrals in \reff{threeint} after partial integration, because there the test function $\vp$ could have compact support.

But now, in  Section \ref{sec4}, i.e. in the case of Robin boundary conditions, the reason for using the cut-off functions $\eta_\ve$ in \reff{barudef1} and \reff{barudef2} is to avoid boundary integrals in \reff{threeint} after partial integration, because now the test function $\vp$ is not allowed to have compact support. Here the cut-off functions $\eta_\ve$ 
are not needed for
certain boundary condition to be satisfied (because the approximate solutions $\bar u_\ve$ are not obliged to satisfy any boundary conditions).

This technical detail is important because of the following: If it would be possible to improve the choices \reff{barudef1} and \reff{barudef2}
of the approximate solutions appropriately (for example, to avoid the use of cut-off functions or to use other smoothing operators), then, perhaps, it would be possible to prove better than \reff{estim1} and \reff{estim2} error estimates.
\end{remark}

\begin{remark}
A similar to Theorem \ref{mainRob} result is true also for elliptic systems with mixed boundary conditions such that each component $u^\al$ of the solution satisfies homogeneous Dirichlet boundary conditions on its own $\al$-depending boundary part (cf. \cite {GrR}).
\end{remark}

\end{document}